\documentclass[10pt,draft]{amsart}

\usepackage{ogonek}

\numberwithin{equation}{section}

\newcommand{\Isom}{\operatorname{Isom}}

\newcommand{\ben}{\begin{enumerate}}
\newcommand{\een}{\end{enumerate}}

\newcommand{\bea}{\begin{eqnarray}}
\newcommand{\ba}{\begin{array}}
\newcommand{\bean}{\begin{eqnarray*}}
\newcommand{\ea}{\end{array}}
\newcommand{\eea}{\end{eqnarray}}
\newcommand{\eean}{\end{eqnarray*}}
\newcommand{\beq}{\begin{equation}}
\newcommand{\eeq}{\end{equation}}
\newcommand{\bthm}{\begin{thm}}
\newcommand{\ethm}{\end{thm}}
\newcommand{\blem}{\begin{lem}}
\newcommand{\elem}{\end{lem}}
\newcommand{\bprop}{\begin{prop}}
\newcommand{\eprop}{\end{prop}}
\newcommand{\bcor}{\begin{cor}}
\newcommand{\ecor}{\end{cor}}
\newcommand{\bdfn}{\begin{dfn}}
\newcommand{\edfn}{\end{dfn}}
\newcommand{\brem}{\begin{rem}}
\newcommand{\erem}{\end{rem}}
\newcommand{\bpf}{\begin{proof}}
\newcommand{\epf}{\end{proof}}
\newcommand{\bfact}{\begin{fact}}
\newcommand{\efact}{\end{fact}}
\newcommand{\bobs}{\begin{obs}}
\newcommand{\eobs}{\end{obs}}
\newcommand{\bexam}{\begin{exam}}
\newcommand{\eexam}{\end{exam}}

\renewcommand{\mod}{\operatorname{mod}}

\alph{enumii} \roman{enumiii}

\newtheorem{thm}{Theorem}[section]
\newtheorem{prop}[thm]{Proposition}
\newtheorem{lem}[thm]{Lemma}

\newtheorem{cor}[thm]{Corollary}
\newtheorem{fact}[thm]{Fact}
\newtheorem{obs}[thm]{Observation}

\theoremstyle{definition}

\newtheorem{dfn}[thm]{Definition}
\newtheorem{rem}[thm]{Remark}
\newtheorem{exam}[thm]{Example}

 \newtheoremstyle{claimstyle}%
   {}%             space above
   {}%             space below
   {\normalfont}%     body font
   {}%                indent
   {\itshape}%        header font
   {.}%               punctuation
   { }%               space after head
   {\thmnote{#3}}%    typeset note only.

\theoremstyle{claimstyle}
\newtheorem*{varclaim}{}

\newenvironment{claim}[1][Claim]{\begin{varclaim}[#1]}{\end{varclaim}}

\newenvironment{subproof}{\begin{proof}}{%
               \end{proof}}
\def\cA{\mathcal A}                    
                    
\def\cL{{\mathcal L}}

\def\N{{\mathbb N}}                \def\Z{{\mathbb Z}}      \def\R{{\mathbb R}}
\def\C{{\mathbb C}}                      \def\oc{{\hat \C}}

\def\Ch{\oc}

                \def\Ga{\Gamma}           
                         
               \def\sg{\sigma}

\newcommand{\eps}{\varepsilon}
\newcommand{\lam}{\lambda}

\newcommand{\ph}{\varphi}
\newcommand{\al}{\alpha}
\newcommand{\ga}{\gamma}

\def\1{1\!\!1}

\def\and{\text{ and }}

\def\Comp{\text{Comp}}

\def\({\bigl(}                \def\){\bigr)}

                        \def\^{\tilde}

\def\cl{\text{cl}}

\def\arg{\text{arg}}

\def\D{{\mathbb D}}
\newcommand{\Hdim}{\operatorname{Hdim}}

\newcommand{\jul}{\mathcal{J}}
\newcommand{\fat}{\mathcal{F}}

\newcommand{\conical}{\Lambda _c}

\newcommand{\rad}{{\mathcal J}_c}

\begin{document}
%***********************************************

\title[Rigidity and absence of line fields]
{Rigidity and absence of line fields \\
 for meromorphic and Ahlfors islands maps}
\date{\today}
% % author information %
\author[\sc Volker MAYER]{\sc Volker MAYER}
%\\[0.5cm]\rm \large Preliminary Version}
\address{Volker Mayer, Universit\'e de Lille I, UFR de Math\'ematiques,
UMR 8524 du CNRS,
59655 Villeneuve d'Ascq Cedex, France}
\email{volker.mayer@math.univ-lille1.fr\newline \hspace*{0.3cm} Web:
math.univ-lille1.fr/$\sim$mayer}
\author[\sc Lasse REMPE]{\sc Lasse REMPE}
%\\[0.5cm]\rm \large Preliminary Version}
\address{Lasse Rempe, Department of Mathematical Sciences, University of Liverpool,
L69 7ZL, United Kingdom}
\email{l.rempe@liverpool.ac.uk\newline \hspace*{0.3cm} Web: http://pcwww.liv.ac.uk/$\sim$lrempe
}
%
% dedication %
%\dedicatory{}
%
% AMS information %
%\thanks{}

\begin{abstract}
 In this note, we give an elementary proof of the absence of 
  invariant line fields on the conical Julia set of
  an analytic function of one variable. This proof applies not only to
  rational as well as transcendental meromorphic functions
  (where it was previously known), but even
  to the extremely general setting of Ahlfors islands maps as defined
  by Adam Epstein. 

 In fact, we prove a more general result on the absence of
  invariant \emph{differentials}, measurable with respect to a 
  conformal measure that is supported on the (unbranched) conical
  Julia set. This
  includes the study of cohomological equations for 
  $\log|f'|$, which are relevant to a number of well-known
  rigidity questions. 
\end{abstract}

\maketitle

%***********************************************
\section{Introduction}

An \emph{invariant line field} of a holomorphic function $f$
 is a measurable field of non-oriented tangent lines that is
 supported on a completely invariant set of positive Lebesgue measure.
 The \emph{no invariant line field conjecture} states that,
 with the exception of flexible Latt\`es mappings (see below),
 rational functions do not support
 invariant line fields on their Julia sets.
 Man\'e, Sad and Sullivan \cite{mss}
 showed that this conjecture is equivalent to the
 density of hyperbolic maps in the 
 space of all rational functions (a conjecture that can be traced back 
 to Fatou; compare \cite[p.\ 59]{mcm}). Correspondingly, the
 problem---which remains open 
 even for quadratic polynomials,
 $f(z)=z^2+c$---has received significant attention. 

However, it is well-known that such line fields do not exist when there
 is sufficient expansion. 
 In particular, the absence of line fields on the
 \emph{conical} (or \emph{radial}) set---on which there is nonuniform 
 hyperbolicity---%%
 was proved by McMullen in \cite{mcm}. This result was generalized to
 arbitrary transcendental meromorphic functions by Rempe and van Strien
 \cite{rvs}, extending earlier partial results by Graczyk, Kotus and
 Swi\k{a}tek \cite{gks}.

These original proofs use orbifold theory and are somewhat
 involved. In \cite{mm}, a rather elementary proof of McMullen's
 result was given (and the argument applies equally to transcendental
 entire functions). In this note, we extend this reasoning to the vast
 class of \emph{Ahlfors islands maps} introduced by Epstein \cite{ep}. These
 are nonconstant
 holomorphic functions $f:W\to X$, where $X$ is a compact
 Riemann surface and $W$ is a nonempty open subset of a compact Riemann surface
 $Y$, that satisfy a certain function-theoretic
 transcendentality condition near every
 point of $\partial W$. (See Section \ref{sec:ahlforsislands} for the precise definition.)
 For the purposes of iteration, one is interested in the case where
 $Y=X$; the set of all Ahlfors islands maps with this property will be
 denoted by $\cA(X)$. Such a function is \emph{non-elementary} if it
 is not a conformal automorphism. 
 We note that this class contains all rational functions $f:\Ch\to\Ch$,
  all transcendental
  meromorphic functions $f:\C\to\Ch$, the functions of \cite{bdh}
   and Epstein's \emph{finite-type maps}. The later class
  includes all (iterated) parabolic renormalizations of rational maps.

 As mentioned above, our results concern the \emph{conical}
   (or \emph{radial}) Julia set. There are several
   inequivalent definitions of this set in the literature.
   For our purposes, we denote by $\conical(f)$ the 
   \emph{(branched) conical set} of $f$,
   that is, the set of points  in the Julia set for which it is possible to 
   pass, using
   iterates of $f$, 
   from arbitrarily small scales to a definite scale by a map of 
   bounded degree $d$. The set of
   points for which we can take $d=1$, i.e. the set of points
   of $\conical(f)$ that have univalent blow ups, is denoted
   $\rad(f)$ and called the \emph{unbranched conical Julia set}. 
   See Section \ref{sec:conical} for the formal definitions.
 
\bthm\label{theo main}
 A non-elementary
  Ahlfors islands map $f\in \cA(X)$ does not support an invariant line field on its
  conical set $\conical(f)$ unless $f$ is either a linear torus endomorphism 
  or a flexible
  Latt\`es map on the Riemann sphere.
\ethm

In particular, our theorem applies to transcendental meromorphic functions,
 and thus provides an elementary
 proof of the result from \cite{rvs}. We also emphasize that
 the argument from \cite{rvs} does not extend to our more general
 setting of Ahlfors islands maps.
As usual, the proof begins by a renormalization argument that shows that
 such a line field would have to be \emph{univalent} near some point of
 the Julia set, and the main step of the proof consists of providing
 a classification of maps that possess such line fields
 (Theorem \ref{thm:univalentlinefield}).

%%% LR: Only small changes up to here.
%%%     In particular notation for the radial Julia set - write
%%%     J_r(f)?
%%%    I don't usually use the curly notation for the Julia set,
%%%     but I don't care particularly. 

\medskip
 
 Similarly as in \cite{my}, our way of establishing Theorem \ref{theo main}
  allows us to handle a much more general framework. 
  For the purpose of this introduction, let us consider the case of
  the unbranched conical set $\rad(f)$, which allows for
  simpler statements. 

 First of all, Theorem \ref{theo main} only applies to maps
  whose Julia set has positive Lebesgue measure. It is natural to replace
 Lebesgue measure by other dynamically significant measures, such as 
 invariant Gibbs states and conformal measures. In the setting of rational 
 functions,  Fisher and Urba\'nski \cite{fu}  studied 
 invariant line fields that are not defined on a set of positive Lebesgue 
  measure 
 but on the support of some conformal measure.
  (More precisely, they work with invariant 
  Gibbs states, but these are equivalent to corresponding conformal measures).
 They obtained  a statement analogous to Theorem \ref{theo main},
 the only difference being that there is a longer list of exceptional
 %%%(one dimensional) 
%%% LR: Why do we specifically mention one-dimensional here? 
 functions. This work has then been used to establish
 ergodicity of the scenery flow of a hyperbolic rational function (\cite{bfu}). 
  
 An equivalent way of stating that a function $f$, defined on a subset of
 the
 complex plane, has an invariant line field is that 
 $\arg\, f'$ is cohomologous to a constant. 
  More precisely, there is a measurable function $\cL$
 such that 
\beq %\label{7}
\cL _{f(z)} = \cL _z + \arg \, f'(z) \quad \mod\ \pi \quad \text{a.e.}
\eeq
 There is another circle of rigidity questions that relies on a different 
  cohomological equation.
  This time, instead of  $\arg\, f'$, the crucial point is to
  study when $\log |f'|$ is cohomologous to a constant.
  This is a central point of the measure rigidity
  problem, initiated by Sullivan \cite{su} and developed 
  further by Przytycki-Urba\'nski \cite{pu}; see also the paper \cite{ku} 
  by Kotus-Urba\'nski 
  where this problem has been investigated for a class of transcendental
  meromorphic functions. 

 Another striking example is Zdunik's paper \cite{zd} which, in combination
  with a result by Makarov \cite{mk} on harmonic measure, shows that 
  Julia sets of polynomials
  are fractals. More precisely, if $f$ is a polynomial with connected Julia 
  set $\jul(f)$
  then either $\jul(f)$ is a circle or line segment 
  (in which case $f(z)=z^{\pm d}$ or $f$ is a Chebyshev polynomial)
  or the Hausdorff dimension $\Hdim (\jul(f) ) >1$.
  The heart of this dichotomy results from stochastic limit theorems
  (central limit theorem and the law of iterated logarithm) and the smooth
  exceptional examples correspond exactly to the case where the variance
  equals zero.
  This last condition is equivalent to the statement
  that a cohomological equation for
  $\log |f'|$ holds $\mu_f$-a.e., where $\mu_f$ is 
  the maximal entropy measure (which, for a polynomial, coincides with 
  harmonic measure of the basin of infinity).

  We give a unified approach to all of these problems.
  In order to do so, and in order to make sense of the questions on
  an arbitrary Riemann surface $X$, we will recast the problem in terms
  of \emph{invariant differentials}. Recall that an
  $(m,n)$-differential, where $m,n\in\Z$ are not both equal to zero, takes
  the form $\mu(z)dz^md\bar{z}^n$ in local coordinates.
  For example, an invariant line field can be identified with an
  invariant Beltrami (i.e., (-1,1)-) differential 
  $\mu\, d\overline{z}/dz$ with $|\mu|=1$ (see \cite[p.47]{mcm}). Also, 
  the statement that 
   $\log |f'|$ is cohomologous to a constant corresponds 
  to saying that some $(1,1)$-differential is invariant 
  \emph{up to a multiplicative constant}. 

 For a more detailed discussion of these concepts, and of 
  conformal measures,
  we refer the reader to Sections  \ref{sec confm} and \ref{sec diffs}.
  Note that the surface $X$ comes equipped with a conformal metric
  of constant curvature: the spherical metric on the sphere, a flat metric
  on a torus, or the hyperbolic metric otherwise, and it is this metric
  that will be used in the notion of a conformal measure.

 %******************************************  ok up to here ******************************

 In this setting, we generalize Theorem \ref{theo main}
  to obtain the following result.
 
 \bthm \label{theo main'}
  Let $f:W_f\to X$ be a non-elementary Ahlfors Island map, with
   $W_f\subset X$. Suppose that, for some $\kappa >0$ and $t\geq 0$,
   there is 
   a 
   $\kappa \|f'\|^t$--conformal probability measure $\nu$ with
   $\nu(\rad(f))=1$. 

 If there exist $(m,n)\in\Z^2\setminus\{0,0\}$ and a $\nu$-measurable
  $(m,n)$-differential that is $f$-invariant
  up to a multiplicative constant, then one of the following
  holds. 

 \begin{itemize}
   \item[(a)] $X=W_f$ is either the Riemann sphere or a torus, and
    $f$ is conformally conjugate to a
     (not necessarily flexible) Latt\`es mapping,
      a linear toral endomorphism, a power map, or
      a Chebyshev polynomial or its negative. 
   \item[(b)] $m+n=0$ and
       there is a finite set $A\subset\jul(f)$ such that
       $\jul(f)$ is locally a $1$-dimensional analytic curve near
       every point of $\jul(f)\setminus A$. 
       In particular,
       $\partial W_f$ is a totally disconnected set. 
      \end{itemize}
 \ethm

 \brem \label{rmk:conformal}
 The condition $\nu (\rad(f)) =1$ 
  is automatically fulfilled for every ergodic invariant measure $\nu$
  of a rational function having strictly positive Lyapunov exponent 
  (see \cite{led}). In particular, this is the
  case for the maximal entropy measure and the equilibrium states mentioned 
  above.

 For a large class of meromorphic 
  functions,  measures supported on
  $\rad(f)$ were constructed in \cite{myurb} that are conformal
  with respect to a suitable conformal metric which is not the spherical
  metric. 
  Theorem \ref{theo main'} also holds for these
  measures (see Definition \ref{defn:compatible} and Lemma
  \ref{lem:compatible}). 
\erem

\brem
An example for case (b) which is not a power map or a
 Chebyshev polynomial is given by $f(z)=z^2-c$ with $c<-2$.
 The Julia set of this polynomial
 is a Cantor set contained in the real line and any constant line field 
 defined on $\jul(f)$
 is $f$--invariant.

 In many situations there is more rigidity in case (b), 
  meaning that one can give the
  precise list of functions that fit into this case 
  (see for example Theorem 1.2 and Claim 3.3 in \cite{my}). 
  Moreover:
\begin{itemize}
\item  If $f$ is rational, then 
 in case (b) the Julia set must be contained
 in a circle on the Riemann sphere \cite{bergweilereremenko,eremenkostrien}. 
\item For entire functions case (b) never happens since the Julia set cannot
   contain isolated Jordan arcs by a result of T\"opfer 
   \cite[Theorem 20]{bergweilersurvey}.
 \item If $\jul(f)$ \emph{is} an analytic curve and $f:\C\to\Ch$ is meromorphic
 with only finitely many critical and asymptotic values, then
 $\jul(f)$ is a straight line \cite{urb}. As far as we know,
 all known examples for case (b) are maps with $X=\Ch$ for which $\jul(f)$ is contained
 in a circle on the Riemann sphere. Examples that are not rational maps are given
 e.g.\ by members of the tangent family. 
 \end{itemize}
\erem

 We also prove Theorem \ref{theo main'} in the case where 
 $f$ has an invariant differential
 that is \emph{continuous} on some relatively open subset of $\jul(f)$.
 (In this case there is no need to assume the existence of
 the measure $\nu$.) In the setting
 of transcendental functions we therefore have the following statements,
 which appear to be new even in the case of invariant line fields.

\begin{cor}
 For an entire function $f$, there is no
  invariant $(m,n)$-differential that 
  is continuous on a relatively open subset of $\jul(f)$. 
  (In particular, $\jul(f)$ does not support
    continuous invariant line fields.) 

 Similarly, if
   a transcendental meromorphic function $f$ has an invariant
  $(m,n)$-diffe\-rential that is continuous on a relatively open
  subset of $\jul(f)$, then $m+n=0$ and 
  $\jul(f)$ is contained in an analytic curve. 
\end{cor}

\section{Ahlfors islands maps}\label{sec:ahlforsislands}

\bdfn \label{defn:ahlfors}
 Let $X$ and $Y$ be compact Riemann surfaces and let $W\subset Y$ be open
 and nonempty. A nonconstant holomorphic function
 $f:W\to X$ is called an \emph{Ahlfors islands map}
  if there is a finite number $k$
  such that the following property holds:

 \emph{Let $V_1,...V_k\subset X$ be Jordan domains with pairwise disjoint closures
 and let
 $U\subset X$ be open and connected with $U\cap \partial W \neq \emptyset$.
 Then, for every component $U_0$ of $U\cap W$, there is $i\in \{1,...,k\}$
 such that $f$ has a simple island over $V_i$ in $U_0$, i.e. there is a
 domain $G\subset U_0$
 such that $f:G\to V_i$ is a conformal isomorphism.}

An Ahlfors islands map is called
 \emph{elementary} if it is a conformal isomorphism.
\edfn

The Ahlfors islands property is implicit in Epstein's thesis
 \cite{adamthesis} and is stated explicitly in \cite{ep}. A more
 detailed study is undertaken in \cite{adamrichard}.
 This class includes all meromorphic functions
 (by the Ahlfors five islands theorem, which gives the class its name)
 as well as the functions considered in \cite{bdh} and the
 \emph{finite-type maps} in the sense of Epstein.
 There are also many examples of Ahlfors islands maps that do not belong
 to these
 categories; compare \cite{exoticbd}.

Given a Riemann surface $X$, we denote the set of all non-elementary
 Ahlfors islands maps $f:W\to X$ with $W\subset X$ by
 $\cA (X)$, and set 
   \[ \cA := \bigcup_X \cA(X). \]
 If $f$ is an Ahlfors islands map, we denote the domain of $f$ by
   $W_f$.

The class of Ahlfors islands maps is closed under composition. (This is easy
 to see when the functions are elementary or their
 domains have nontrivial boundaries. The only non-elementary Ahlfors islands
 maps without boundary are rational maps and endomorphisms of the torus, and these
 can be treated classically.)
 In particular,
 if $f\in\cA(X)$, then the iterates $f^k$ belong to $\cA(X)$ for all
 $k\geq 1$.

 The Fatou set $\fat(f)$ of a map $f\in\cA$ is defined as the set of 
  points $z\in X$ that have
  a neighborhood $U$ with the property that either all iterates are defined on
  $U$ and
  form a normal family there or $f^n(U)\subset X\setminus\cl(W_f)$ for some 
  $n\geq 0$.
  The Julia set $\jul(f)$ is the complement of $\fat(f)$ in $X$.
  All results from the basic iteration theory of transcendental entire and meromorphic functions
  also hold for Ahlfors islands maps. In particular, we have the following
  important result (proved by Epstein \cite{ep} along the lines of Baker's original proof
  \cite{baker} for transcendental entire functions).

\bthm[Baker, Epstein] \label{repellers dense}
 Repelling periodic points are dense in the Julia set of a non-elementary
  Ahlfors islands map $f\in\cA$ and the Julia set $\jul(f)$ is a perfect subset of $X$.
\ethm

% For further details, compare \cite{re}.
%*************************** next lemma needed for random only ********************************
The key fact in the proof of the density of repellers is the following lemma
(see \cite[Lemma 2.5]{re}).

\blem \label{Islands near jul}
Let $f\in \cA (X)$. Then there is a number $k$ with the following property: if
$V_1,...,V_k \subset X$ are Jordan domains with disjoint closures, then there is
some $V_j$ such that every open set $U$ with $U\cap \jul(f) \neq \emptyset$ contains
an island of $f^n$ over $V_j$, for some $n$.
\elem

%*************************** end lemma needed for random only ********************************

Suppose that $X$ is a compact Riemann surface, $W\subset X$ is open and
 $f:W\to X$ is analytic. If $z_0\in W$ is a repelling fixed point of $f$, then
 there is a map $\Psi$, defined and univalent in a neighborhood of $0$ and satisfying
 $\Psi(0)=z_0$, such that
   \begin{equation}  \label{eqn:linearizing}
    f(\Psi(z)) = \Psi(\lambda z), \end{equation}
 where
 $\lambda = f'(z_0)$. (The map $\Psi$ is unique up to precomposition by a linear map.)

 Using \eqref{eqn:linearizing}, the map $\Psi$ can be extended
  to some maximal domain $W_{\Psi}\subset\C$. We refer to this maximal extension
  $\Psi:W_{\Psi}\to\C$ as a \emph{Poincar\'e function} associated to $z_0$ and $f$. The
  following was first observed by Epstein.

\begin{lem}\label{11}
 If $f$ is an Ahlfors islands map, then
  $\Psi$ is also an Ahlfors islands map.
\end{lem}
\begin{proof}
 Let $V$ be a neighborhood of $0$ on which $\Psi$ is univalent. By definition,
  a point $z$ belongs to $\partial W_{\Psi}$ if and only if there is $n$ such that
  $z/\lambda^n \in V$ and $\Psi(z/\lambda^n)\in \partial W_{f^n}$. Since the Ahlfors islands
  property is preserved under iteration,
  and since $\Psi$ is univalent on $V$, the claim follows.
\end{proof}

One can also define \emph{Picard points maps} or \emph{Casorati-Weierstra{\ss} maps} in analogy
 to the definition of Ahlfors islands maps, generalizing these classical theorems of complex
 analysis instead of the Ahlfors islands theorem; compare \cite{adamrichard}. However, these
 do not seem to be sufficiently strong to obtain an interesting dynamical theory.
 For our purposes, the following definition will nonetheless be useful.

\bdfn
 Let $X$ and $Y$ be compact Riemann surfaces, let 
  $W\subset Y$ be open and nonempty and
  $f:W\to X$ be holomorphic and non-constant. We say that $f$ has the
  \emph{weak Casorati-Weierstra{\ss} property} if, for every open set 
  $U\subset Y$ with
  $U\cap \partial W\neq \emptyset$, the image $f(U\cap W)$ is dense in $X$.
\edfn

Clearly, if $f$ is an Ahlfors islands map and $U$ is a component of $W_f$, then the restriction
 $f_{|U}$ is also an Ahlfors islands map, and in particular a weak Casorati-Weierstra{\ss} map.
We note furthermore that Lemma \ref{11} holds also for the class of weak 
 Casorati-Weierstra{\ss} mappings.
 The reason is, again, that the
 weak Casorati-Weierstra{\ss} property is preserved under iteration.

\section{Conical set and renormalization}\label{sec:conical}

In the following, we always suppose that $X$ is a compact Riemann surface and 
 that
 $f\in \cA (X)$ is a Ahlfors islands map. We also fix a metric of 
 constant curvature
 on $X$---the spherical metric if $X=\Ch$, a flat metric (say of area $1$) 
 if $X$ is a torus,
 and the hyperbolic metric otherwise.
 A disk of radius $r>0$ and centered at a point $z\in X$ (with respect to
 this metric) will be denoted by 
 $D(z,r)$. (If $r$ is sufficiently small, depending on $X$, the set
 $D(z,r)$ is simply 
 connected, and hence conformally equivalent to the standard unit disk.)
 Euclidean disks will be denoted
 by $\D (z,r) \subset \C$. If $z\in\Omega\subset X$, the connected
 component of $\Omega$ containing $z$ will be denoted by
 $\Comp_z (\Omega )$.

 The literature contains various different definitions of the 
  \emph{conical Julia set}. 
  (The term ``radial'' is also used interchangeably with ``conical''.) 
  In many cases, the
  dynamical difference between the different definitions 
  is not significant (the reader is referred to 
   the discussion given in \cite{P} and also in \cite{rvs}). 
  We shall use both the \emph{unbranched} conical set $\rad(f)$ and the
  \emph{branched} conical set $\conical(f)$. We note that
  $\conical(f)$ is the most general definition that appears
  in the literature.
    
\bdfn \label{conical}
 The $d$-branched conical Julia set $\rad^{(d)}(f)$ is the set of
  points $z\in \jul(f)$ for which 
 there exists $\delta >0$ and $n_j\to\infty$ such that 
 $D_j=\Comp _z \left(f^{-{n_j}}(D(f^{n_j}(z),\delta )) \right)$ is simply 
 connected and
\beq \label{zoom} f^{n_j} : D_j \longrightarrow D(f^{n_j}(z), \delta )\eeq
 is defined and a proper map of degree at most $d$.

 The \emph{unbranched conical set} is denoted by
  $\rad(f) := \rad^{(1)}(f)$. The \emph{branched conical set} of $f$ is
  defined as
   \[ \conical(f):=\bigcup_{d\geq 1} \rad ^{(d)}. \]
\edfn

Every repelling periodic point belongs to $\rad(f)$.  
 In particular, $\rad$ and $\conical$ are
 dense (and forward-invariant) subsets of $\jul(f)$.
 The main feature of conical points is that one can make good 
 renormalizations, passing from small to large scales.
 Moreover, these "good" renormalizations characterize conical points 
 \cite[Proposition 2.3]{h}:

\bfact \label{lemm renormalisation}
 Let $z_0\in\jul(f)$, and let $\phi:U\to \C$ be a
  local chart in a neighborhood of $z_0$, with $\phi(z_0)=0$.

 Then 
  $z_0\in\conical(f)$ if and only if there 
  are integers $n_j\to\infty$ and a sequence $r_j>0$ with 
  $r_j\to 0$ such that 
\beq\label{8}
\Psi _j (z)= f^{n_j}(\phi^{-1}(r_j\cdot z))
\eeq
converges uniformly on the unit disk $\D \subset \C$ to a non-constant holomorphic map 
 \[ \Psi : \D \longrightarrow \Omega=\Psi (U )\subset X. \]

 Moreover, $z_0$ belongs to $\rad(f)$ if and only if
  $0$ is not a critical point of $\Psi$.
\efact

\section{Conformal measures} \label{sec confm}
 Conformal measures were introduced by Sullivan (in analogy with the
  case of Kleinian groups) as natural substitutes of 
  Lebesgue measure in the case
  where the Julia set has zero area. 

\bdfn
Let $f\in \cA (X)$ and let $\kappa >0$, $t\geq0$. Let
  $m$ be a probability measure supported on 
  on $\jul(f)$.

Then $m$ is called \emph{$\kappa \|f'\|^t$-conformal}
 if, for every Borel set $E\in \jul(f)$
 for which $f_{|E}$ is injective, 
 \[ m(f(E))=\int _E \kappa \|f'\|^t \, dm \;\;. \]
 (Here the derivative is measured with respect to the natural
  conformal metric on $X$.)
\edfn

%Usually conformal measures are required to be
% probability measures, but here we will allow the measures
% to be infinite. This is natural since
% the Julia set of an Ahlfors islands map is not necessarily compact.

The number $\log \kappa$ is most often the topological pressure. 
 There are two important examples
 in the case of a rational function $f$. First of all, if $t=0$ and 
 $\kappa= deg(f)$, then the (invariant) maximal entropy
 measure $\mu_f$ is $deg(f)|f'|^0=deg(f)$-conformal. 
 The second interesting case is when $\kappa =1$.
 Sullivan showed that there is always a $|f'|^t$--conformal measure for
 a minimal exponent $t>0$ \cite[Theorem 3]{su2}.
%%% LR: Ich glaube, Sullivan hat schon die Existenz eines
%%%      minimalen Exponenten gezeigt?

\subsection*{Conformal measures and the conical set}
 Conformal measures are particularly useful when they are supported on
  the unbranched conical set $\rad(f)$: here we can pass from
  small scales to large scales using univalent iterates, and 
  conformality means that the measure behaves in a well-controlled
  manner under these blow-ups. We now formalize 
  the key property we require.

 Recall that a value $w\in X$ is \emph{Fatou-exceptional} for
  $f\in\cA(X)$ if the set of iterated pre-images of
  $w$ is finite; we denote the set of Fatou-exceptional values by
  $E_F(f)$. For example, $\infty$ is a Fatou-exceptional value for
  any transcendental entire function $f$. Let us also denote
  by $E_B(f)$ the set of \emph{branch-exceptional values} for $f$
  in the sense of \cite{rvs}; i.e.\ the set of points that have 
  only finitely many \emph{unbranched} preimages. For example, 
  $0$ and $\infty$ are branch-exceptional points for 
  $f(z)=z^2$. The Ahlfors islands property
  implies that $E_B(f)$ and $E_F(f)$ are finite. 

 \bdfn\label{defn:compatible}
   Let $f\in\cA(X)$. A measure $m$ is called
    \emph{$\rad(f)$-compatible} if the following hold.
  \begin{enumerate}
   \item[(a)]
 The topological support of $m$ contains the Julia set
      $\jul(f)$, and furthermore $m(\{z\})=0$ for all $z$.
   \item[(b)] $m$ is locally finite, except possibly near points of
     the branch-exceptional set $E_B(f)$.
   \item[(c)] $m(\rad(f))>0$.
   \item[(d)] Let $w\in \C\setminus E_B(f)$, 
     suppose that $\delta>0$ is such that 
     $D := D(w,\delta)$ has finite $m$-measure and set
     $\tilde{D} := D(w,\delta/2)$. Then for every
     $\eta>0$, there exist $\eps>0$ and $n_0\in\N$ such that the following
     holds. 

   Let
     $U$ be a component of $f^{-n}(D)$, for some $n\geq n_0$, such that
     $f^n:U\to D$ is univalent, and let 
     $\tilde{U}\subset U$ be a component of $f^{-n}(\tilde{D})$.
     Then for any subset $A\subset \tilde{U}$ with
       \[ m(A)/m(\tilde{U})<\eps, \]
     we have 
     \[ m(f^n(A))/m(\tilde{D}) < \eta. \]
  \end{enumerate}

 Similarly, we say that $m$ is 
  \emph{$\conical(f)$-compatible} if it satisfies condition
   (a) and furthermore:
  \begin{enumerate}
   \item[(b')] $m$ is locally finite except possibly near points of
   the Fatou-exceptional set $E_F(f)$,  
   \item[(c')] $m(\conical(f))>0$, 
   \item[(d')] for every $\Delta>0$, condition (d) holds with 
      ``$w\in\C\setminus E_B(f)$'' replaced by ``$w\in\C\setminus E_F(f)$'' and 
      ``$f^n:U\to D$ is univalent'' 
      replaced by ``$f^n:U\to D$ is proper of degree at most $\Delta$''. 
  \end{enumerate}
 \edfn

\begin{lem} \label{lem:compatible}
  Let $\mu$ be a  $\kappa \|f'\|^t$-conformal measure with 
   $\mu(\rad(f))>0$ and $\mu(\{x\})=0$ for all $x\in X$.
   Then $\mu$ is $\rad(f)$-compatible.

  The conformal measures constructed in \cite{myurb} are 
   $\rad(f)$-compatible.
  
  Hausdorff measure of any dimension (including Lebesgue measure)
   is $\conical(f)$-compatible provided that
   $\conical(f)$ has positive and locally finite Hausdorff measure of
   the corresponding dimension. 
\end{lem}

\begin{proof}
Concerning property (d),  it follows for $\kappa \|f'\|^t$-conformal measures
directly from the Koebe distortion theorem. For Hausdorff measures property (d')
can be shown by contradiction and using a renormalization argument like
the one of the proof of Lemma \ref{lem univalent}.

The manuscript \cite{myurb} concerns a quite general class of hyperbolic meromorphic functions
and provides, in particular, conformal probability measures associated to H\"older continuous potentials 
and a particular choice of Riemannian metric. 
Lemma 5.20 of \cite{myurb} yields property (d) and implies that these measures do not have
any mass on points. Furthermore, 
Proposition 5.21 of \cite{myurb} shows that the support of these mesures 
is the set of points that do not escape to infinity which,
the functions being (topologically) hyperbolic, 
is a subset of the unbranched conical set $\rad (f)$.
This shows property (c).
\end{proof} 

\begin{lem}
  If $m$ is a $\rad(f)$- or $\conical(f)$-compatible measure, then
   $m$ is ergodic,  and
   $\rad(f)$ resp.\ $\conical(f)$ has full $m$-measure
   $m$-almost every point as a dense orbit in $J(f)$.

  (In particular, if $\conical(f)$ has positive Lebesgue measure, then
   $\jul(f)=X$ and $f$ is ergodic.) 
\end{lem}
\begin{proof}
 The proof is analogous to the proof for Lebesgue measure; see 
  \cite{mcm, rvs}. Since the result is stated in a somewhat
  unfamiliar framework, let us give some of the details.

 Let $F$ be a forward-invariant subset of $\rad(f)$ with $m(F)>0$.
  Because $m$ is locally a finite Borel measure, we can pick a density point $z_0$ of
  $F$. That is, the density of $F$ within small simply-connected sets of 
  bounded
  geometry around $z_0$ is
  close to $1$. 

 Since $z_0$ belongs to $\rad(f)$, there is a disk $D=D(w,\delta)$ 
  and univalent pull-backs
  $U_j$ of $D$, $f^{n_j}(U_n)=D$, with $n_j\to\infty$ and such that 
  $f^{n_j}(z_0)\to w$. We may assume that $z_0$ is not a repelling periodic
  point that belongs to the branch exceptional set, in which case $w$ also
  does not belong to the branch-exceptional set. So 
  we can choose $\delta$  sufficiently small that
  $m(D)<\infty$. 

 Set $\tilde{U}_j := \Comp_{z_0}(f^{-n_j}(\tilde{D}))$, where 
  $\tilde{D}=D(w,\delta/2)$.
  Then, by choice of $z_0$, we have
    \[ m(\tilde{U}_j\cap F)/m(\tilde{U}_j)\to 1, \]
  and hence, by the definition of a $\rad(f)$-compatible measure, and since
   $F$ is forward-invariant, 
    \[ m(\tilde{D}\cap F)/m(\tilde{D}) \to 1. \]
   So $F$ has full measure in $\tilde{D}$. 

 We claim that it follows that $m(X\setminus F)=0$. Indeed, otherwise we can
  (again using the density point theorem, and the fact that $m$
   does not give positive measure to singletons) find an infinite, pairwise
   disjoint collection of balls in each of 
   which $F$ does not have full measure.
   Using the Ahlfors islands property (or, more precisely, 
   Lemma \ref{Islands near jul}),
   at least one of these balls must have an iterated univalent 
   preimage in $D$. But then $F$ has full measure inside this iterated 
   preimage, again by the definition of a $\rad$-compatible measure. 
   This is a contradiction. 

 So in particular $\rad(f)$ has full $m$-measure and $m$ is ergodic.
  Furthermore, for any open set $U$ intersecting the Julia set, the
  set of points whose orbits do not enter $U$ is forward-invariant,
  and hence has measure zero by the above. It follows that $m$-almost
  every orbit is dense in $J(f)$. 
\end{proof} 

In \cite{fu}, Fisher and Urba\'nski
 worked with invariant Gibbs states of H\"older continuous potentials $\Phi$ 
 satisfying
 the condition $\sup \Phi > P(\Phi)$, the pressure of $\Phi$. It is well known 
 that the H\"older 
 continuity of $\Phi$ combined with exponential shrinking of inverse branches 
 yields distortion
 estimates. The exponential shrinking of inverse branches is the object of 
 Man\'e's theorem.
 We would like to mention that Man\'e's theorem is also available for certain 
  meromorphic functions
 \cite{rvs}. Therefore we could also extend our work to such measures.

\section{Invariant and univalent differentials}\label{sec diffs}

 Let $(m,n)\in\Z^2\setminus \{(0,0)\}$. A $(m,n)$-differential $\mu$ on the
  Riemann surface $X$ is 
  a differential that is locally of the form
  $\mu(z) dz^m d\bar{z}^n$. More precisely, the differential
  is defined by a collection of functions
  $\mu_i:U_i\to\C$, where $(U_i)$ is a collection of
  local charts with corresponding local coordinates
  $z_i=\ph_i(p)$, such that
\beq\label{1.1}\mu_i = \mu_j \left(\frac{dz_j}{dz_i}\right)^m
              \left(\frac{d\overline{z}_j}{d\overline{z}_i}\right)^n 
     \quad  \text{holds on} \quad U_i\cap U_j\;.
 \eeq
 All $(m,n)$-differentials we consider are assumed 
  to be non-zero by convention.
  We will be interested in differentials that are continuous
  on an open and dense subset of the Julia set, or those that are 
  measurable with respect to a given (conformal) measure. The pullback
  $f^*\mu$
  of a differential under a holomorphic function is well-defined except
  at critical points. We shall be interested in differentials that 
  are \emph{invariant up to a multiplicative constant}, by which we mean
  that there is a constant $c\in\C$ such that
\beq \label{1.2} f^*\mu =c\mu \eeq
  holds where defined (for continuous differentials) resp.\
  almost everywhere (for measurable differentials). 
  Note that the definition implies that the support of $\mu$ is backward
  invariant up to the countable set of critical points resp.\ 
  up to a set of measure zero. 

\bexam As explained in \cite[p.47]{mcm}, a line field can be identified with a unit Beltrami differential, i.e. a
$(-1,1)$--differential whose functions $\mu_i$ have modulus one. A line field is invariant
if the corresponding Beltrami differential is invariant. 
 Examples of functions that have an invariant measurable line field (with respect to Lebesgue measure)
 are given by (flexible) Latt\`es functions and torus endomorphisms,
 power maps and Chebychev polynomials.
\eexam

\bexam
 All Latt\`es maps (even rigid ones) have differentials that are Lebesgue a.e.-invariant up to
  a multiplicative constant. For example, consider the case of a map that is the quotient of
  a linear map under the full symmetry group of a hexagonal lattice. The differential
  $dz^3/d\bar{z}^3$ is invariant under this symmetry group, and hence descends to a 
  $(3,-3)$-differential (defined and continuous except at finitely many points) that is
  invariant up to a multiplicative constant under the quotient map. 
\eexam

\bexam
Suppose that a $(1,1)$--differential $\mu$ is $f$--invariant up to multiplication, where
  $f$ is analytic in the complex plane. Then the invariance equation becomes
  $\mu\circ f |f'|^2=c\mu$ (valid in some local charts).
  Taking the logarithm of this equation gives exactly the condition that $\log|f'|$ is cohomologous
  to a constant. 
\eexam

 In \cite{mcm} (and \cite{rvs}), the notion of a locally univalent line field is used. This
  concept has the following straightforward adaption to differentials.
  We say that a $(m,n)$--differential $\mu$ is \emph{univalent} near a point $z$ if there exists a
  neighborhood $V$ of $z$ and a
  conformal map $\psi:U\to V$ such that the pull back of 
  $\mu$ by $\psi$ is a constant differential, i.e.
$$ \psi ^* \mu= c\,dz^m\,d\overline{z}^n \quad on \quad U\,.$$
  Similarly, we shall say that a differential $\mu$ is locally univalent 
  on a set $A$ if every $z\in A$ has a neighborhood $U$ such that
  $\mu|_{A\cap U}$ agrees with a univalent differential. 

 We leave the proofs of the following simple fact to the reader.

\bobs \label{obs:constantdifferential}
  Let $(m,n)\in\Z^2\setminus\{(0,0)\}$. Suppose that $f$ is a holomorphic 
   function such that
   the constant differential
   $dz^md\bar{z}^n$ is invariant up to a multiplicative constant. 
   Then $f$ is affine:
   $f(z)=az+b$ for some $a,b\in\C$, $a\neq 0$. 
\eobs

% \begin{lem}
%  Let $(m,n)\in\Z^2$ with $m+n\neq 0$. Let $\lambda\in\R$ with $|\lambda|>1$ and let 
%   $\mu$ be a locally univalent $(m,n)$-differential. If the restriction of
%   $\mu$ to the real axis is invariant, up to
%   a multiplicative constant, under $z\mapsto \lambda z$, then $\mu$ is constant. 
% \end{lem}

The first step towards the rigidity theorems is a well known fact that renormalization at 
a conical point at which the differential is continuous in measure leads to locally univalence
of the differential. 

\blem \label{lem univalent}
 Let $f\in\cA (X)$ and let $(m,n)\in\Z^2\setminus\{(0,0)\}$.
 Suppose that $\mu$ is a $(m,n)$-differential supported on the Julia set
 such that either
 \ben 
 \item $\mu$ is invariant up to a multiplicative constant and 
   $\mu|_{\jul(f)}$ is 
   continuous at some point of $\conical(f)$, or
 \item $\mu$ is measurable with respect to some $\rad(f)$-compatible or
   $\conical(f)$-compatible measure $\nu$, and $\mu$ is 
   invariant up to a multiplicative constant $\nu$-a.e.
 \een
  Then $\mu$ is univalent on a nonempty and relatively open
  subset of $\jul(f)$.
\elem
\begin{proof}
This is a standard fact, proved by
 renormalizing at a conical point $z_0$ at which the differential
  is continuous (resp.\ continuous in measure). So, let $z_0$ be such a point and let
 \beq \label{14}
\Psi_j =f^{n_j}\circ \al _j \to \Psi : \D \to \Omega = \Psi (\D ) \subset X
\eeq
 be the renormalization given by Fact \ref{lemm renormalisation}.
  From the assumption that $\mu$ is continuous (in measure) at $z_0$, and the
  definition of $\rad(f)$-\ resp.\ $\conical(f)$-compatible measures,
   it follows 
 that the pull back
$\psi^* \mu$ is a constant differential 
$\Psi ^{-1} (\conical \cap \Omega )$ ($\nu$-a.e.).
(Compare 
\cite[p. 4359]{mm} and
\cite[p. 561]{my}.)

This shows that $\mu$ is univalent at $z_0$ when $z_0\in \rad$.
Finally, if $z_0\in \conical$ and if $\Psi$ is not univalent, then we can conclude as follows.
Since $\Omega \cap \jul(f)\neq \emptyset$ and since $\jul(f)$ is a perfect set
there is a conformal restriction of $\Psi$ whose image still intersects $\jul(f)$.
This shows that the differential is univalent near some point in the Julia set.
\end{proof}

\section{The case of Lebesque measure: Proof of Theorem \ref{theo main}}

 In this section, we study the case where the function supports a 
  differential that is invariant (up to a multiplicative constant) and
  univalent on a full complex neighborhood of a point of $\jul(f)$. Indeed,
  the following result completely classifies the set of such functions.

\bthm \label{thm:univalentlinefield}
  Let $(m,n)\in\Z^2\setminus\{(0,0)\}$ and 
   suppose that $f\in \cA (X)$ has an $(m,n)$-differential 
    (not necessarily supported on the Julia set)
   that is invariant
   up to a multiplicative constant and univalent on a 
   full complex neighborhood of some point of
   $\jul(f)$. 

 Then $X=W_f$ is either the Riemann sphere or a torus, and $f$ is conformally 
  conjugate to one of the following:
 \begin{enumerate}
  \item A power map $z\mapsto z^j$, $j\in\Z$, $|j|\geq 2$;
  \item a Chebyshev polynomial or its negative;
  \item a (rigid or flexible) Latt\`es mapping;
  \item a linear toral endomorphism.
 \end{enumerate}
\ethm

 By Lemma \ref{lem univalent}, this result implies Theorem \ref{theo main},
  and indeed more generally it implies Theorem \ref{theo main'} when 
  $\jul(f)=X$. (The remaining case is considered in 
  the next section.)  

To begin the proof, suppose that $f$ satisfies the hypotheses of the Theorem
  \ref{thm:univalentlinefield} and let 
  $\mu$ be the invariant differential. By assumption, 
  there is some open connected set $U$, intersecting the Julia set,
  and a univalent map $\Psi:V\to U$ such that
  $\Psi^*(\mu)=dz^m\,d\overline{z}^n$.
 
The set $U$ contains some repelling periodic point $p$ of $f$; we may assume 
 (by postcomposing
 $\Psi$ with a translation) that $\Psi^{-1}(p)=0$.
 Let us pass to an iterate $g=f^k$ for which $p$ is fixed. We recall that again 
 $g\in \cA(X)$;
 note that $g$ also preserves $\mu$. Consequently, the function
\beq\label{sigma} \sigma = \Psi ^{-1} \circ g \circ \Psi,\eeq
defined locally near the origin, preserves the constant differential $dz^md\overline{z}^n$
up to multiplication. Thus, by Observation
  \ref{obs:constantdifferential}, we have 
  $\sigma (z)= \lam z $, where $\lam =g'(p)$ has modulus greater than $1$.
  Hence $\Psi$ linearizes $g$ near $p$, 
  and we can extend it using the functional relation (\ref{eqn:linearizing}) 
  to a maximal domain $W_{\Psi}\subset\C$. By the functional relation,
  the pullback of $\mu$ under this extended map is a constant differential.

%%%?? In the following we identify $\sigma $ and $\lam$.

So far, this is exactly the same argument as in \cite{mm}.
For rational or even entire functions the map $\Psi$ is defined globally on $\C$, and
 in \cite{mm} it is then shown directly that $\Psi$ is an exponential, trigonometric or
 elliptic function.
For meromorphic non-entire and thus for general Ahlfors islands maps, it is no longer
a priori clear that $W_{\Psi}=\C$.
This is the point where we need new arguments in order to carry over the ideas from \cite{mm, my}.
We start with the following lemma.

\blem\label{10}
 Let $X$ be a Riemann surface and $V\subset\C$ be connected, open and nonempty.
  Suppose that $\Psi:V\to X$ is a holomorphic map such that the pushforward of 
  $dz^md\overline{z}^n$ is (almost everywhere) well-defined.
  Then
 \begin{enumerate}
   \item[(DT)] for all $z_1,z_2\in V$ with $\Psi(z_1)=\Psi(z_2)$, there exists a
    ``deck transformation''
    $$\gamma(z) = \al z + \beta\quad , \quad \al\in\C\setminus\{0\}\; ,\;\; \beta\in\C\; ,$$
     such that
    $\gamma (z_1)=z_2$ and $\Psi \circ \gamma =\Psi$.
 \end{enumerate}
\elem

\brem \label{rema 1}
In the case of a line field, i.e. if $m=1$ and $n=-1$, 
 an easy calculation show that the numbers $\al$ for the maps in the deck 
 transformation property (DT) are real. This has some importance in 
 Proposition \ref{prop:orbifold}
 and in the proof of Theorem \ref{theo main}.
\erem

\begin{proof}
 Let us first suppose that $z_2$ is not a critical point of $\Psi$. 
  Then we can locally
  define an inverse branch $\Psi_2 ^{-1}$ mapping
  $\Psi (z_2)$ to $z_2$. Hence we can define near $z_1$ a map
  $\ga = \Psi_2 ^{-1} \circ \Psi$.
  Clearly $\ga (z_1)=z_2$. The assumption means that $\ga$ preserves the 
  constant differential,
  and hence
  it is of the desired form. Because $V$ is connected, the relation
  $\Psi\circ\gamma=\Psi$ holds on all of $V$ by the identity theorem.

 If $z_2$ is a critical point, then we can perturb $z_1$ and $z_2$ to nearby
  points $z_1'$ and $z_2'$ with $\Psi(z_1')=\Psi(z_2')$, and obtain a map $\gamma$ with
  $\gamma(z_1')=z_2'$ as above.
  Provided the perturbation was small enough, it follows that also $\gamma(z_1)=z_2$.
\end{proof}

Now we have essentially reduced the problem to the following proposition, which has a distinctly
 classical flavor.

\begin{prop} \label{prop:orbifold}
 Suppose that $V\subset\C$ is nonempty, open, and connected, that $X$ is a compact Riemann
  surface and that
  $\Psi:V\to X$ is a holomorphic map satisfying the weak
  Casorati-Weierstra{\ss} property and
  the deck-transformation property (DT) as above.

 Then $V=\C$, $X$ is either the Riemann sphere or a torus, and $\Psi$ is one of the following,
  up to pre- and postcomposition by conformal isomorphisms:
  \begin{enumerate}
   \item the exponential map $\exp:\C\to\C\setminus\{0\}$;
   \item the cosine map $\cos:\C\to\C$;
   \item a projection $\pi:\C\to X$, where $X$ is a torus;
   \item a Weierstra{\ss} $\wp$-function $\wp:\C\to \Ch$.
   \item the map $\wp'$ or $(\wp')^2$, where $\wp$ is the Weierstra{\ss} function
     associated to a hexagonal lattice;
   \item the map $\wp^2$, where $\wp$ is the Weierstra{\ss} function associated to
     a square lattice. 
  \end{enumerate}
 Furthermore, if $\alpha\in\R$ for the maps in the property (DT), then the last two cases
  do not occur. 
\end{prop}
\begin{proof}
Let $\Gamma$ be the set of all mappings
 $\gamma(z) = \al z + \beta$, $\al\in\C\setminus\{0\}$, $\beta\in\C$ satisfying
  $\Psi\circ\gamma=\Psi$ and $\gamma(V)\cap V\neq\emptyset$.
  Because of (DT) and the weak Casorati-Weierstra{\ss} property,
 $\Gamma$ is quite a rich set.

\begin{claim}[Claim 1]
The domain $V$ is $\Gamma$--stable. That is, $\gamma(V)=V$ for every
  $\gamma\in \Gamma$.
\end{claim}
\begin{subproof}
Let $\ga \in \Gamma$.
The set $U=\ga (V) \cap V$ is an open and non-empty subset of $V$
on which we have $\Psi = \Psi \circ \ga ^{-1}$. It follows that $\Psi$ (and $\Psi \circ \ga ^{-1}$)
have an analytic extension to the connected open set $V\cup \gamma (V )$.
Hence  $V\cup \gamma (V )\subset V$ since otherwise we have a contradiction
to the fact that $\Psi$ does not extend beyond $V$ by the
weak Casorati-Weierstra{\ss} property.
\end{subproof}
\begin{claim}[Claim 2]
 $\Gamma$ is a discrete subgroup of $\Isom(\C)$.
\end{claim}
\begin{subproof}
The fact that $V$ is invariant under the elements of $\Ga$ immediately implies that
 $\Ga$ is closed under composition, and the inverse of any element of $\Gamma$ again
 belongs to $\Gamma$ by definition. Hence $\Gamma$ is a subgroup of the group of non-constant
 affine maps $\C\to\C$.

This group is discrete because otherwise there would be a sequence
 $\ga _j\in \Ga $ such that $\ga_j \to Id$ and this is in contradiction with the fact that
 $\Psi ^{-1} (z)$ is a discrete subset of $V$ for every $z\in\C$.

 We showed that $\Ga$ is a discrete group, hence it must be a group of isometries.
\end{subproof}

\begin{claim}[Claim 3]
 We have $V=\C$.
\end{claim}
\begin{subproof}
Suppose that the domain $V$ had some finite boundary point $a\in \partial V$. Let
 $z_0\in V$ and let $\delta$ be such that $\D(z_0,2\delta )\subset V$. By
 the weak Casorati-Weierstra{\ss} property, there is a point $z_2\in \D(a,\delta)$ such that
 $\Psi(z_2)\in \Psi(\D(z_0,\delta))$; so $\Psi(z_2)=\Psi(z_1)$ for some $z_1\in \D(z_0,\delta)$.
 By property (DT), there is $\gamma\in\Gamma$ with
 $\gamma(z_1)=z_2$, and as we have just shown, $\gamma$ is an isometry. But then
  \[ a\in \D(z_2,\delta) = \gamma(\D(z_1,\delta))\subset \gamma(\D(z_0,2\delta)) \subset \gamma(V), \]
 which contradicts Claim 1.
\end{subproof}

 Since $V=\C$, $\Gamma$ is a group of isometries and $\Psi$ has infinite degree,
  we have the following possibilities,
  up to a change of coordinates (according to the classification of wallpaper groups
  without reflections):
  \begin{enumerate}
   \item $\Gamma=\{z+2\pi i k; k\in\Z\}$ and $\Psi=\exp$;
   \item $\Gamma=\{\pm z + 2\pi k; k\in\Z\}$ and $\Psi = \cos$;
   \item $\Gamma$ is a lattice and $\Psi$ is the projection to the corresponding torus;
   \item $\Gamma$ is the product of a lattice and the order two subgroup
     $\{z\mapsto\pm z\}$, and $\Psi$ is a Weierstra{\ss} $\wp$-function.
   \item $\Gamma$ is either the full symmetry group or a rank-two subgroup of
     a hexagonal lattice, and $\Psi=\wp'$ or $\Psi=(\wp')^2$, where $\wp$ is the
     corresponding Weierstra{\ss} function.
   \item $\Gamma$ is the symmetry group of a square lattice, and $\Psi=\wp^2$, where
    $\wp$ is the corresponding Weierstra{\ss} function.
  \end{enumerate}
 In the latter two cases, the group $\Gamma$ contains rotations of order greater
  than two, and hence these cannot occur if all deck transformations are
  of the form $\alpha z + \beta$ with $\alpha\in\R\setminus\{0\}$. 
\end{proof}

\begin{proof}[Conclusion of the proof of Theorem \ref{thm:univalentlinefield}]
 Let $\Psi$ be the Poincar\'e function as above, and $V=W_{\Psi}$. Then it follows
  that $V=W_{\Psi}=\C$, and that $\Psi$ is one of the maps in the finite list above. Since
  $\Psi$ semi-conjugates a linear map to $g$, this proves the theorem for the iterate
  $g$ of $f$.

 For the details that the same is true of $f$,
  we again refer to \cite{mm}
  since, as soon as we get $V=\C$, we are back to the situation of a global map $\Psi$.
  In fact, one uses the global map $\psi$ to lift a convenably chosen restriction of $f$.
  This lift again preserves a constant differential up to multiplication, 
  from which follows that the lift is affine.
  It follows that $\Psi$ semi-conjugates this affine map and $f$ which implies that 
  $f$ 
  is of the desired form.
\end{proof}

\begin{proof}[Proof of Theorem \ref{theo main}]
Suppose that $f$ has an invariant line field a.e. on a set of positive 
 Lebesgue measure.
 By Lemma \ref{lem univalent}, this line field is univalent in a neighborhood of 
 some point
 in the Julia set of $f$. So it follows from 
 Theorem \ref{thm:univalentlinefield} combined with 
 Remark \ref{rema 1} that
 we either have the desired conclusion, or $f$ is a power map or a Chebyshev 
 polynomial. However,
 the latter maps have Julia sets of zero measure (a circle or a line segment, 
 respectively),
 so we are done.
\end{proof}

%************************************************************************************************
\subsection*{Further remarks}
 We note that  we have not used in the proof of Theorem \ref{theo main}
 the full strength of the Ahlfors islands property. What really is needed is the 
 weak Casorati-Weierstra{\ss} property
 combined with the density of repellers in the argument just before \eqref{sigma}.
 The latter is not known for Casorati-Weierstra{\ss} functions. The result is 
 nevertheless true in this more general setting since we now 
  explain how  the density of repellers can be avoided in that proof. 
 
 \bthm \label{thm:CW}
 Theorem \ref{theo main} is true for every function
  $f:W_f\to X $ that has the weak Casorati-Weierstra{\ss} property.
  Similarly, Theorem \ref{thm:univalentlinefield} holds for such functions
  if ``some point of $\jul(f)$'' is replaced by
  ``some point of $\conical(f)$''. 
\ethm
 
\begin{proof}
We have to inspect the situation given in \eqref{14} more closely.
So let again $\Psi=\lim_{j\to \infty}f^{n_j}\circ \al_j:\D\to \Omega$ be a limit
of the renormalization at the conical point $z_0\in \conical$, where
$z_0$ is chosen to be a Lebesgue continuity point (for the case of Theorem
 \ref{theo main}) or a point near which the differential $\mu$ is
 univalent. In both cases, it follows as in Lemma \ref{lem univalent} that
 the pullback of $\mu$ under $\Psi$ is constant. 
 It suffices to show 
 that $f$ has a repelling periodic point in $\Omega$. 

We can choose $V\subset \Omega$ to be simply connected such that $\Psi (0)\in V$ 
 and such that
 $\Psi : U\to V$ is a proper map of finite degree say $d$, where 
 $U$ is the connected
 component of $\Psi ^{-1} (V)$ that contains $0$.

Given $j,k$, we denote $g= f^{n_{j+k}-n_j}$. 
 If $j,k$ are big enough then $g^{-1} (V)$ has a component $V_0$ that is relatively
compact in $V$ and that contains $f^{n_j}(z_0)=f^{n_j}\circ \al_j (0)$. 
 The map $g:V_0\to V$ is proper.

Let $U_0$ be a connected component of $\Psi^{-1}(V_0)$; then 
 $\overline{U}_0\subset U$.
 Because the set of critical values of $\Psi$ is countable, we can pick
 $a\in U_0$ such that $g(\Psi(a))$ is not a critical value of $\Psi$;
 let $b$ be a preimage of $g(\Psi(a))$. 
 We can define again locally a function $\sg$
 such that $\Psi \circ \sg =g\circ \Psi$ and $\sg (a)=b$.
Considering once more the invariant differential up to multiplication it follows that 
$\sg (z)=\al z+\beta$ is a globally defined affine map. We can therefore set $U_1=\sg (U_0)$.

\

\noindent
{\bf Claim:} $U_1=U$.

\

Let us admit the Claim for a moment. Then we have that $U_0$ is relatively compact in $U_1$
and it follows that the affine map $\sg : U_0\to U_1$ has a repelling fixed point in $U_0$.
Then $g$ must also have a repelling fixed point in $V_0\subset \Omega$ and we are done.

It remains to establish the Claim. 
Indeed, if $z\in U_0$ such that $\sg (z)\in \partial U$ then we immediately have a contradiction
with $\Psi\circ \sg (z)= g\circ \Psi (z)\in V$. Therefore we have that $U_1\cap \partial U=\emptyset$.
On the other hand, $U_1\cap U\neq \emptyset$ and so $U_1\subset U$.
Moreover $\Psi =g\circ \Psi \circ \sg ^{-1} :U_1\to V$ is a proper map and so 
 $\partial U_1\subset \partial U$. In conclusion we have inclusion of the domains and 
 of their boundaries. They must be equal.

\end{proof}
%**********************************************************************************************

\section{The case of conformal measures: proof of Theorem \ref{theo main'}}
\label{sec:cohom}

To prove Theorem \ref{theo main'}, we still need to deal with the case
 where the Julia set might not equal the entire Riemann surface $X$. 
 In this case, Lemma \ref{lem univalent} still shows that the restriction
 of the differential $\mu$ to the Julia set is locally univalent. 
 If the Julia set is not locally contained in an analytic curve, then
 it follows that the differential $\mu$ can be extended to a locally univalent
 invariant differential also outside the Julia set, and the result from the
 previous section applies. Otherwise, it is possible to see that
 we must have $m+n=0$. We now make these arguments somewhat more precise. 

\bthm \label{theo 7.1}
  Let $(m,n)\in\Z^2\setminus\{(0,0)\}$ and 
   suppose that $f\in \cA (X)$ has an $(m,n)$-differential $\mu$, supported
   on the Julia set, that is invariant
   up to a multiplicative constant and univalent on a 
   relative neighborhood of some point of
   $\jul(f)$. 

  Then either 
  \begin{enumerate}
   \item $f$ is conformally conjugate to a power map, a Chebyshev
     polynomial or its negative, a Latt\`es mapping or a 
     toral endomorphism, or \label{item:lattescase}
   \item 
     $m+n=0$ and
       there is a finite set $A\subset\jul(f)$ such 
       $\jul(f)$ is locally a $1$-dimensional analytic curve near
       every point of $\jul(f)\setminus A$. In particular,
       $\partial W_f$ is a totally disconnected set. 
  \end{enumerate}
\ethm
\begin{proof}
 Similarly as in the previous section, let $z_0$ be a periodic point
  for $f$ near which the restriction of the
  differential $\mu$ to the Julia set is univalent. We may also suppose
  for simplicity that $z_0$ is not branch-exceptional, and 
  pass to an iterate $g$ for which $z_0$ is fixed. We prove the claim
  for $g$; the conclusion for $f$ then follows as in the previous
  section. Set $\lambda := g'(z_0)$ and let
  $\Psi:W_\Psi\to X$ be the Poincar\'e function at $z_0$. Also set
  $\hat{\jul}=\Psi^{-1}(\jul(f))$ and let $\nu$ be the univalent
  $(m,n)$-differential
  on $\hat{\jul}$ obtained from $\mu$ by pullback under $\Psi$. 
  Note that  $\nu$ is univalent near zero by assumption, and hence
  univalent on all of $\hat \jul$ by the functional relation
  \eqref{eqn:linearizing} and the fact that $\mu$ is invariant up to
  a multiplicative constant.

 The fact that $\nu$ is invariant under multiplication by 
  $\lambda$, up to a multiplicative constant $C$, implies that
  $\nu$ is a constant differential. Indeed, we have
    \[ \lambda^{m}\bar{\lambda}^n \cdot \nu(z) = C\cdot \nu(\lambda z) \]
  on $\hat{\jul}$.
  First substituting $z=0$, we see that 
   $\nu(z)=\nu(\lambda z)$ for all $z\in\hat{\jul}$ and, by the
   identity theorem, $\nu$ is locally constant near $0$. Finally,
   the functional relation again implies that $\nu$ is constant everywhere.
   So without loss of generality, we can suppose that 
   $\nu = dz^m\,d\bar{z}^n$.
   
   Since the support of $\mu$
    is contained in the Julia set we consider the following modified deck 
    transformation condition:
 \begin{enumerate}
   \item[(DT')] for all $z_1,z_2\in \hat{\jul}$ with $\Psi(z_1)=\Psi(z_2)$, there exists a
    ``deck transformation''
    $$\gamma(z) = \al z + \beta\quad , \quad \al\in\C\setminus\{0\}\; ,\;\; \beta\in\C\; ,$$
     such that
    $\gamma (z_1)=z_2$ and $\Psi \circ \gamma =\Psi$.
 \end{enumerate}

\medskip

We begin by noting that we are done whenever
  (DT') holds:

 \begin{claim}[Claim ]
 If (DT') holds, then the formally stronger condition (DT) is also satisfied. 
  In particular (by Proposition \ref{prop:orbifold}) 
   we are in Case (\ref{item:lattescase}) of Theorem \ref{theo 7.1}.
\end{claim}
\begin{subproof}
 We inspect the proof of Proposition  \ref{prop:orbifold}, letting
  $\Gamma$ be the group of all affine deck transformations as in
  (DT').
  The arguments proving 
  Claims 1 and 2 apply without any changes and yield
  that the domain $W_\Psi $ is $\Gamma$--stable and that $\Gamma$ is a 
  discrete subgroup
  of $\Isom(\C )$. 
  Similarly, the Ahlfors islands property, and the fact
  that $\hat\jul$ is a perfect set (compare Theorem \ref{repellers dense}
  and Lemma \ref{Islands near jul})
  easily imply that there is a point in $\jul(f)$ that has 
  $\Psi$-preimages
  near every point of $W_{\Psi}$. This implies that Claim 3
  also holds; i.e.\ $W_{\Psi}=\C$. Hence $\Gamma$ is one of the
  finitely many groups described at the end of the Theorem.

 Since $\Psi$ semiconjugates $f$ and $z\mapsto \lambda z$, and
  $\hat{\jul}$ is completely invariant under multiplication by $\lambda$,
  it follows
  that $\lambda\Gamma\subset\Gamma$. Thus the restriction
  $f|_{\jul(f)}$ has finite degree $|\lambda|$ or $|\lambda|^2$. 
  By the Ahlfors islands property, it follows that
  $\partial W_{f}=\emptyset$---i.e., $f$ is either a rational map
  or a toral endomorphism---and $f$ itself has finite degree.

 Now we can verify property (DT). Take points $z_1$ and
  $z_2$ with $w:=\Psi(z_1)=\Psi(z_2)$; as in Lemma \ref{10}, we may assume
  that these are not critical points for $\Psi$. Let
  $\psi$ be the branch of $\Psi^{-1}\circ\Psi$ that takes $z_1$ to
  $z_2$; we must show that $\psi$ is affine. 
  If $\Psi(z_1)\in \jul(f)$, then we are done by (DT'). So
  suppose that $w$ belongs to a component $U$ of
  $\fat(F)$. If $U$ is not a rotation
  domain, then $\mathcal{P}(f)\cap U$ is either discrete in $U$ or
  has at most one accumulation point,
  where $\mathcal{P}(f)$ is the postcritical set of $f$
  (i.e., the closure of the set of critical orbits). Hence 
  we can choose a curve $\gamma$ connecting
  $\Psi(z_1)$ to a point $\tilde{w}\in \jul(f)$ such that
  $\gamma$ does not contain any postcritical points of $f$.
  Then any branch of $\Psi^{-1}$ can be analytically extended
  along $\gamma$, and hence we obtain a branch 
  $\tilde{\psi}$ of $\Psi^{-1}\circ\Psi$ that is defined on a neighborhood
  of a point in the Julia set. By (DT'), the map $\tilde{\psi}$, and hence
  $\psi$, is affine.

 Finally, suppose that $U$ is a rotation domain. For
  $i\in\{1,2\}$ and $k\in\N$, define
  $\zeta_i^k := \Psi(z_i/\lambda^k)$. Let $k$ be sufficiently large that
  $\zeta_1^k, \zeta_2^k\notin U$, and let $\gamma$ be a curve 
  connecting $w$ to $\jul(f)$ such that $\gamma$ is disjoint
  from $f^k(\mathcal{P}(f)\setminus U)$ and the set of
  critical values of $f^k$. Recalling that
  $\Psi(z)=f^k(\Psi(z/\lambda^k))$, we see again that the
  branches of $\Psi^{-1}$ taking $w$ to $z_1$ and $z_2$ can
  be continued analytically along $\gamma$, and we obtain the
  conclusion as above.

 (Alternatively, we could have observed that 
   $\hat{\jul}$ is the limit of $\Gamma(0)/\lambda^n$ as
   $n\to\infty$, and that hence $\jul(f)$ must contain an analytic curve.
   Furthermore, by the Gross star theorem \cite[Page 292]{nevanlinna},
   every branch of $\Psi^{-1}$ can be continued analytically
   along almost every radial line. Together, these facts lead to a shorter
   but 
   less dynamical proof of
   property (DT).)
\end{subproof}
  
 In order to check whether (DT') holds or not, let $z_1,z_2\in \hat{\jul}$ and suppose
 that $\gamma$ is a local defined holomorphic function with $\gamma (z_1)=z_2$
 and such that $\Psi\circ\gamma = \Psi$. Such a map preserves
  the constant differential $\nu$ when restricted to $\hat{\jul}$. If it would
  preserve the differential everywhere (not just on $\hat{\jul}$), then
  the map is affine and we are done. Here however we have to distinguish a number of cases.
  
  \smallskip

 \noindent
 {\textbf{$\jul(f)$ is not locally contained in an analytic curve near $z_0$.}}
  Since $z_0$ was assumed not to be branch-exceptional, it follows that
  $\jul(f)$ cannot locally be an analytic curve anywhere, and hence the
  same is true for $\hat{\jul}$. Any deck transformation $\gamma$
  as above satisfies
    \begin{equation} \label{eqn:decktransinvariant}
      \gamma'(z)^m\,\cdot\, \bar{\gamma'(z)}^n = 1 \end{equation}
  for all $z\in\hat{\jul}$. By assumption, $\hat{\jul}$ is not
  locally contained in an analytic curve, and hence 
  $\gamma'$ must be constant; i.e.\ $\gamma$ is affine, as desired.

 \smallskip

 \noindent
 {\textbf{$\jul(f)$ is locally contained in an analytic curve near $z_0$, and
   $m+n\neq 0$.}}
  Then $\hat{\jul}$ is an analytic curve near $0$; since $\hat{\jul}$ is
   invariant under multiplication by $\lambda$, it follows that
   $\lambda$ is real and that 
   $\hat{\jul}$ is contained in a straight line through the origin;
   without loss of generality, we may assume that $\hat{\jul}\subset\R$. 

  In particular, every deck transformation $\gamma$ as above must
   preserve the real line; i.e.\ $\gamma'(z)\in\R$ for all $z\in J(f)$. 
   Given that $m+n\neq 0$, we see that
   \eqref{eqn:decktransinvariant} implies that
   $\gamma'(z)^{m+n}$, and hence $\gamma'(z)$, is constant.

 \smallskip
 {\textbf{$\jul(f)$ is locally contained in an analytic curve near $z_0$, and
   $m+n=0$.}}
  Again, we can assume that $\hat{\jul}\subset\R$. Recall that
   $\Psi$ is an Ahlfors islands map. Hence, if $z\in J(f)$, then
   $J(f)$ is locally an analytic curve near $z$ provided that
   there is $w$ with $\Psi(w)=z$ and $\Psi'(w)\neq 0$. 

  Because $\Psi$ is an Ahlfors islands map, such a point $w$ exists
   for all $z\in J(f)$ apart possibly 
   from a finite set $A$. Furthermore, recall that the boundary
   $\partial W_f$ of the domain of definition of $f$ is contained
   in the Julia set of $f$. If $\partial W_f$ contained a nontrivial
   subcontinuum, then there would be a point where
   $\partial W_f$ is an analytic curve that is isolated in the Julia set,
   which is impossible by the Ahlfors islands property. Hence 
   $\partial W_f$ is totally disconnected.
\end{proof}

%*************************************************************************************************

%********************************************************************************************************

%\end{document}

\end{document}